\documentstyle[12pt]{article}
\newlength{\abstractwidth}
\renewcommand{\thefootnote}{\fnsymbol{footnote}}
\renewcommand{\thanks}[1]{\footnote{#1}} 
\newcommand{\starttext}{
\setcounter{footnote}{0}
\renewcommand{\thefootnote}{\arabic{footnote}}}

\newcommand{\be}{\begin{equation}}
\newcommand{\bea}{\begin{eqnarray}}
\newcommand{\eea}{\end{eqnarray}}
\newcommand{\ee}{\end{equation}}

\def\ba{\begin{eqnarray}}
\def\ea{\end{eqnarray}}

\def\log{\,{\rm log}\,}

\def\ddb{{\partial\bar\partial}}

\def\[{{\bf [}}
\def\]{{\bf ]}}

\begin{document}
\starttext
\baselineskip=18pt
\setcounter{footnote}{0}

\begin{center}
{\Large \bf On Multiplier Hermitian Structures on Compact
K\"ahler Manifolds}

\bigskip\bigskip

{\large Qi Li}  \\

\bigskip

 Department of Mathematics\\
Columbia University, New York, NY 10027\\

\end{center}

\baselineskip=15pt
\setcounter{equation}{0}
\setcounter{footnote}{0}

\section{Introduction}
\setcounter{equation}{0}

In \cite{M4}, Mabuchi introduced the notion of a  multiplier
Hermitian structure on K\"{a}hler manifolds and a generalization
of the notions of K\"{a}hler-Einstein metric and K\"ahler-Ricci
soliton. In this note we study these new notions on compact
K\"{a}hler manifolds whose first Chern class is positive.

\medskip

Let $M$ be a compact K\"{a}hler manifold of complex dimension $n$
with positive first Chern class $c_{1}(M)>0$ and let $\mathcal{K}$
denote the set of all K\"{a}hler forms $\omega$ on $M$ in the
class $c_{1}(M)$. Assume that $X$ is a holomorphic vector field
on $M$ and that
\begin{eqnarray*}
{\mathcal{K}}_{X}:=\{\omega \in {\mathcal{K}}:
L_{X_{\mathbf{R}}}\omega=0\}\neq \emptyset
\end{eqnarray*}
where ${X_{\mathbf{R}}}=X+\overline{X}$ denotes the real vector
field on $M$ associated to $X$. In this note, all K\"{a}hler
metrics considered will be in this set.  We assume also that $X$
is Hamiltonian, i.e., that we can find a function ${u_{\omega}}\in
C^{\infty}(M)_{\mathbf{R}}$  normalized by
$\int_{M}{u_{\omega}\omega^{n}}=0$ so that
\be
X^\alpha=\frac{1}{\sqrt{-1}}g^{\alpha\overline{\beta}}\partial_{\bar
\beta}u_{\omega}
\ee
where
$\omega=\sqrt{-1}\sum_{\alpha,\beta}g_{\overline{\beta}\alpha}dz^{\alpha}\wedge
d\bar z^{\beta}$. In \cite{FM}, Futaki and Mabuchi proved
that
$$
l_{0}:=\min_{M}u_{\omega},\ l_{1}:=\max_{M}u_{\omega}
$$
are
independent of the choice of $\omega \in {\mathcal{K}}_{X}$.

\medskip

Let $\sigma$ be a real-valued smooth function defined on
interval $[l_{0},l_{1}]$ satisfying one of the following
conditions:
\medskip

(a)  $\dot{\sigma}\leq 0 \leq \ddot{\sigma}$

(b)  $\ddot{\sigma}>0$

\medskip

Here $\dot{\sigma}$ and $\ddot{\sigma}$ are the first derivative
and second derivative of $\sigma$.  Associated to this $\sigma$
and $X$, Mabuchi introduced the following
generalization of the notions of
K\"ahler-Einstein metric and K\"ahler-Ricci soliton,

\medskip
\noindent
{\bf Definition 1}. {\it Let $M$ be a K\"ahler manifold with $c_1(M)>0$
and a holomorphic vector field $X$. Fix a real-valued function $\sigma$ as above.
A metric ${\omega}$ in the class $c_1(M)$
is said to be an Einstein-Mabuchi metric of
type $X$ and $\sigma$ if}
\be
Ric({\omega})+\sqrt{-1}{\partial}\bar{\partial}\sigma(u_{\omega})=\omega.
\ee

\bigskip
\noindent
${\mathbf{Remark}}$: In the definition,
$Ric({\omega})+\sqrt{-1}{\partial}\bar{\partial}\sigma(u_{\omega})$
can be viewed as the Ricci curvature of the metric
${\rm exp}(-\frac{\sigma(u_\omega)}{n})\omega$ which is the multiplier
Hermitian metric introduced in [8].

\medskip
Special cases of Einstein-Mabuchi metrics include:

\medskip

(1) K\"ahler-Einstein metrics, corresponding to $\sigma=0$;

(2) K\"ahler-Ricci soliton, defined by $Ric(\omega)-\omega=L_{X}\omega$,
where $L_{X}$ is the Lie derivative along $X$. This corresponds to
the Einstein-Mabuchi metric of
type $\sigma(s)=-s+C$

(3) Let $h_\omega$ be the Ricci potential, defined by
$Ric(\omega)-\omega=\sqrt{-1}{\partial}\bar{\partial}h_{\omega}$.
If $1-e^{h_{\omega}}$ defines a holomorphic vector field as in
(1.1), then the metric $\omega$ is called the generalized K\"ahler-Einstein
metric with nonvanishing Futaki invariant \cite{M1}.
This metric corresponds to the Einstein-Mabuchi metric of type $\sigma(s)=-\log(s+C)$
where $C$ is a constant strictly greater than
$l_{0}$.

\bigskip

According to a well-known conjecture of Yau \cite{Yau}, the existence of
K\"ahler-Einstein metrics should be equivalent to a notion of stability
in geometric invariant theory. Formulations of versions of this conjecture in
terms of the notion of $K$-stability have been given by Tian \cite{T}
and Donaldson \cite{Donaldson}.
Analytically, the existence of K\"ahler-Einstein
metrics is related to the properness of the functional
$F_{\omega}$. By properness we mean that for any sequence
$\phi_i \in C^{\infty}(M)_{\mathbf{R}}$ such that
$\omega_{\phi_i}=\omega+\sqrt{-1}\partial\overline{\partial}\phi_i>0$,
we must have
$ \limsup_{i\rightarrow \infty}F_{\omega}(\phi_{i})=+\infty$
whenever ${\lim_{i\rightarrow
\infty}}J_{\omega}(\phi_{i})=+\infty$. Here
\begin{eqnarray*}
J_{\omega}(\phi)&:
=&\frac{1}{V} \int_{0}^{1}\int_{M}\dot{\phi_{s}}({\omega}^n-{\omega}_{\phi_{s}}^n)ds\\
F_{\omega}(\phi)&: =&
J_{\omega}(\phi)-\frac{1}{V}\int_{M}\phi{\omega}^n-\log(\frac{1}{V}\int_{M}e^{{h}_{\omega}-\phi}{\omega}^n)
\end{eqnarray*}
for all function $\phi$ such that $\omega_{\phi}>0$. In the
definition of $J_{\omega}(\phi)$, $\phi_s$ is a path connecting
$0$ and $\phi$ with $\phi_0=0$ and $\phi_1=\phi$. We shall also
require the functional $I_{\omega}(\phi)$, which is closely related to
$J_{\omega}(\phi)$ and is defined by
\begin{eqnarray*}
I_{\omega}(\phi)&:
=&\frac{1}{V}\int_{M}\phi({\omega}^n-{\omega}_{\phi}^n).
\end{eqnarray*}
\medskip
In \cite{T}, Tian proved that there exists a K\"ahler-Einstein
metric $\omega_{KE}$ on the K\"ahler manifold $(M, \omega)$ with
$c_{1}(M)>0$ as long as the functional
$F_{\omega}$ is proper. More precisely,  he proved an inequality
of Moser-Trudinger type for K\"ahler-Einstein manifolds $M$ without
nontrivial holomorphic vector fields, i.e., \be
{F}_{\omega_{KE}}(\phi)\geq A {J}_{\omega_{KE}}(\phi)^\gamma-B \ee
where $\gamma={\frac{e^{-n}}{8n+8+e^{-n}}}$. Clearly, this inequality
implies that the  functional ${F}_{\omega_{KE}}(\phi)$ is proper.
Recently it was proved in \cite{PSSW}
by Phong-Song-Sturm-Weinkove that the exponent $\gamma$ can be
taken to be 1. In \cite{CTZ},  the results of \cite{T} have been extended,
under some additional assumptions, to
the case of K\"ahler-Ricci solitons by Cao-Tian-Zhu with
$\gamma=1/4n+5$ for the generalized functionals $\tilde{I}$,
$\tilde{J}$, $\tilde{F}$ associated to the vector field $X$.

\medskip

In this note we generalize the results in \cite{T}, \cite{CTZ},
\cite{PSSW} and \cite{TZ1} to Einstein-Mabuchi metrics. First, we
introduce the appropriate generalizations of the functionals $I$,
$J$ and $F$, which we still denote by $\tilde{I}_{\omega}(\phi),
\tilde{J}_{\omega}(\phi), \tilde{F}_{\omega}(\phi)$. Note that
$\tilde{I}_{\omega}(\phi)$ and $\tilde{J}_{\omega}(\phi)$ first
appeared in \cite{M4}. A key feature of these generalizations is
the use of the volume form $ e^{-\sigma(u_\omega)}\omega^n$
instead of the volume form $\omega^n$. Without loss of generality,
we may assume that
$V=\int_{M}e^{-\sigma(u_\omega)}{\omega}^n=\int_{M}{\omega}^n$. We
set \be
\tilde{I}_{\omega}(\phi):=\frac{1}{V}\int_{M}\phi(e^{-\sigma(u_\omega)}{\omega}^n-e^{-\sigma(u_{\omega_{\phi}})}{\omega}_{\phi}^n)
\ee \be
\tilde{J}_{\omega}(\phi):=\frac{1}{V}\int_{0}^{1}\int_{M}\dot{\phi_{s}}(
e^{-\sigma(u_\omega)}{\omega}^n-e^{-\sigma(u_{\omega_{\phi_s}})}{\omega}_{\phi_s}^n)ds
\ee\be
\tilde{F}_{\omega}(\phi):=\tilde{J}_{\omega}(\phi)-\frac{1}{V}\int_{M}\phi
e^{-\sigma(u_\omega)}{\omega}^n-\log(\frac{1}{V}\int_{M}e^{{h}_{\omega}-\phi}{\omega}^n)\ee
The variational derivative of the functional
$\tilde{F}_{\omega}(\phi)$ is readily computed \be {\delta{\tilde
{F}}}=\int_M
\delta\phi(e^{-\sigma(u_{\omega_{\phi}})}{\omega}_{\phi}^n-
\frac{1}{\int_{M}e^{{h}_{\omega}-\phi}{\omega}^n}e^{{h}_{\omega}-\phi}{\omega}^n)\ee
Thus the critical points $\phi$ of the functional
$\tilde{F}_{\omega}(\phi)$ are given by the equation \be
e^{-\sigma(u_{\omega_{\phi}})}{\omega}_{\phi}^n-
\frac{V}{\int_{M}e^{{h}_{\omega}-\phi}{\omega}^n}e^{{h}_{\omega}-\phi}{\omega}^n=0
\ee This is the equation for Einstein-Mabuchi metrics as we shall
see in the next section.

\bigskip

Let $Aut^0(M)$ be the
identity component of the group of all holomorphic automorphisms
of $M$, and let $G\subset Aut^0(M)$ be a maximal compact subgroup.
Let $Z(X)$ be the compact subgroup of $G$ consisting of
all $g\in G$ such that $Ad(g)X=X$, and let $Z^0(X)$ be the identity
component of $Z(X)$. Let ${\mathcal{H}}_{X}$ be also the set of
all ${X_{\mathbf{R}}}$ invariant functions $\phi$ in
$C^{\infty}(M)_{{\mathbf{R}}}$ such that $\omega_{\phi}$ is in
${\mathcal{K}}_{X}$. In \cite{M4}, Mabuchi proved that
Einstein-Mabuchi metrics on $M$ with respect to $X$ must be
$Z^0(X)$-invariant.  We introduce the
following definition of properness.

\bigskip
\noindent
{\bf Definition 2} {\it The functional
$\tilde{F}_{\omega}$ is said to be proper
with respect to the functional $\tilde J_\omega$ if for any sequence
$\{\phi_i\}$ of $Z^0(X)$-invariant functions with
$\omega_{\phi_i}\in {\mathcal{K}}_X$,
we have
$\limsup_{i\rightarrow
\infty}\tilde{F}_{\omega}(\phi_{i})=+\infty$ whenever
${\lim_{i\rightarrow
\infty}}\tilde{J}_{\omega}(\phi_{i})=+\infty$.}

\bigskip
In this note, we will establish the following
theorems.

\bigskip
\noindent
{\bf Theorem 1} {\it If the functional $\tilde{F}_{\omega}(\phi)$
is proper with respect to the functional $\tilde{J}_{\omega}(\phi)$,
then there exists an Einstein-Mabuchi metric on the K\"ahler
manifold $(M, \omega)$.}

\bigskip
\noindent
{\bf Theorem 2} {\it Let $M$ be a compact K\"ahler manifold with
holomorphic vector field X which admits a Einstein-Mabuchi metric
${\omega}_{EM}$ of type $\sigma$. Assume that $K\subseteq
{Z^0(X)}$ is a closed subgroup whose centralizer in G is finite,
then there are two positive constants A and B such that for any
K-invariant function $\phi$ in ${\mathcal{H}}_{X}$, \be \tilde
{F}_{\omega_{EM}}(\phi)\geq A \tilde {J}_{\omega_{EM}}(\phi)-B
 \ee}

\bigskip
\noindent
{\bf Remark:} The condition that $K\subseteq {Z^0(X)}$ is a closed
subgroup whose centralizer in $G$ is finite is a natural
generalization of a condition introduced in
\cite{PSSW} for the case of K\"ahler-Einstein manifolds with
nontrivial holomorphic vector fields.

\bigskip
\noindent The organization of the note is as follows. In section 2
we review some basic properties of multiplier Hermitian structures
and prove Theorem 1. In section 3 we prove Theorem 2 following the
method of \cite{T}. In details, we follow closely the exposition
of \cite{PSSW}. In the last section we construct a holomorphic
invariant of Futaki type which can be viewed as an obstruction to
the existence of the Einstein-Mabuchi metric.

\section{Proof of Theorem 1}
\setcounter{equation}{0}

To an arbitrary smooth path
${\varphi=\{\phi_{t}; 0\leq t \leq 1}\}$ in ${\mathcal{H}}_{X}$,
it corresponds to  a one-parameter family of K\"ahler forms
$\omega(t)$ in ${\mathcal{K}}_{X}$ by \be
\omega(t):=\omega_{\phi_{t}}=\omega+\sqrt{-1}{\partial}\bar{\partial}\phi_{t},
\ \ 0\leq t \leq 1\ee for $\omega \in {\mathcal{K}}_{X}$.

\medskip

In \cite{M4}, Mabuchi used the method of continuity to deform a
given metric along a path to find the Einstein-Mabuchi metric,
i.e., \be
Ric(\omega_{\phi_t})+\sqrt{-1}{\partial}\bar{\partial}\sigma(u_{\omega_{\phi_t}})=(1-t)\omega+t\omega_{\phi_{t}}
 \ee
which is equivalent to \be
-\sqrt{-1}{\partial}\bar{\partial}\log({\omega}_{\phi_t}^n)+\sqrt{-1}{\partial}\bar{\partial}\sigma(u_{\omega_{\phi_t}})=\sqrt{-1}{\partial}\bar{\partial}t\phi_t
-\sqrt{-1}{\partial}\bar{\partial}\log({\omega}^n)-\sqrt{-1}{\partial}\bar{\partial}{h}_{\omega}\ee
. Then we get the following  complex Monge-Amp\`{e}re equation ,
\be
\frac{{\omega}_{\phi_{t}}^n}{{\omega}^n}=e^{{h}_{\omega}-t\phi+\sigma(u_{\omega_{\phi_t}})}\ee
We want to solve the equation at t=1 which gives the
Einstein-Mabuchi
 metric. Since the right hand side of (2.4) is uniformly bounded when $t=0$, one can solve the equation
at $t=0$ by standard argument for complex Monge-Amp\`{e}re
equation. The implicit function theorem implies that the set
$T=\{t:$ there is a solution at t, $0\leq t \leq 1\}$ is open. If
one can get a uniform $C^{0}$ estimate for $\phi_{t}$ for $t\in
T$, by the standard argument for complex Monge-Amp\`{e}re equation
one can show the closeness of the set. As in \cite{A}, the $C^{0}$
estimate is closely related to the functional $\tilde{I}$,
$\tilde{J}$, $\tilde{F}$. In \cite{M4}, Mabuchi has proved the following results,\\\\
{\bf Theorem [M1]:} {\it For $t\in[1/2,1]$, we have positive real
constants $C_{0}$, $C_{1}$, independent of the choice of the pair
$(\omega_{\phi_t}, t)$ such that\be osc(\phi_{t})\leq
C_{0}(\tilde{I}_{\omega}-\tilde{J}_{\omega})(\phi_{t})+C_1\ee for
all $(\omega_{\phi_t}, t)$ with $t\in[1/2,1]$. Here
$osc(\phi_{t})=\max(\phi_t) -\min(\phi_t)$.}

\medskip

So to prove Theorem 1, we only need to prove that the properness
of the functional $\tilde{F}$ implies that one can derive upper
bound for functional $\tilde{I}-\tilde{J}$. Before proving the
theorem, we review some basic properties of the
multiplier Hermitian structure.

\medskip

Since we consider the volume form
$e^{-\sigma(u_\omega)}{\omega}^n$ in the functional, we need the
following formula for integration by parts, \be
-\int_{M}(\bar{\partial}u,
\bar{\partial}v)_{\omega}e^{-\sigma(u_\omega)}{\omega}^n=\int_{M}(u\overline{(\Delta_{\omega}+\sqrt{-1}\dot{\sigma}(u_{\omega})\bar{X})v})
e^{-\sigma(u_\omega)}{\omega}^n \ \ \  \ee for any complex-valued
smooth function $u$, $v$ on $M$, where
$\Delta_{\omega}=\sum_{\alpha,\beta}g^{\alpha\overline{\beta}}{\partial_\alpha}\bar{\partial}_{\beta}$
is the Lapalacian
 operator for $\omega$.
For simplicity we will use the operator $\tilde{\Box}_{\omega}$ to
denote $\Delta_{\omega}+\sqrt{-1}\dot{\sigma}(u_{\omega})\bar{X}$.

\medskip

One has $u_{\omega_{t}}=u_{\omega}+\sqrt{-1}X(\phi_{t})$. By using
the fact that $\phi_{t}$ is ${X_{\mathbf{R}}}$ invariant, we
have\be u_{\omega_{t}}=u_{\omega}-\sqrt{-1}\bar{X}\phi_{t} \ee Let
$\dot{\phi_{t}}$ denote the partial derivative of $\phi(t)$ with
respect to t, by using the (2.7), it is easy to verify that, \be
\frac{\partial}{\partial{t}}(e^{-\sigma(u_\omega)}{\omega}^n)=
(\tilde{\Box}_{\omega_t}\dot{\phi_{t}})e^{-\sigma(u_\omega)}{\omega}^n\ee
\be
\int_{M}e^{-\sigma(u_{\omega_{\phi_t}})}{\omega}_{t}^n=V=\int_{M}e^{-\sigma(u_\omega)}{\omega}^n
\ \ \ for\ all\ \omega\in {\mathcal{H}}_{X}
 \ee

 Mabuchi also proved the following properties of
the generalized I, J functional:

\medskip \noindent {\bf Fact 1}: {\it$0\leq
\tilde{I}_{\omega}(\phi_{t}) \leq
(m+2)(\tilde{I}_{\omega}-\tilde{J}_{\omega})(\phi_{t})) \leq
(m+1)\tilde{I}_{\omega}(\phi_{t})$, where m is a constant
depending only on $\sigma. $}\\
{\bf Fact 2}: {\it Along the equation, one has\be
\frac{d}{dt}(\tilde{I}_{\omega_{0}}-\tilde{J}_{\omega_{0}})(\phi_{t}))=-\int_{M}(\phi_{t}\tilde{\Box}_{\omega_{\phi_{t}}}\dot{\phi_{t}})e^{-\sigma(u_{\omega_{\phi_t}})}{\omega}_{t}^n
=\int_{M}\{\dot{\phi_{t}}+\tilde{\Box}_{\omega_{\phi_{t}}}\dot{\phi_{t}}\}(\tilde{\Box}_{\omega_{\phi_{t}}}\dot{\phi_{t}})e^{-\sigma(u_{\omega_{\phi_t}})}{\omega}_{t}^n\geq
0\ee i.e., $(\tilde{I}_{\omega}-\tilde{J}_{\omega})(\phi_{t})$ is
increasing along the equation}.

\medskip

We refer interested readers to \cite{M4} for details of the above
properties. We also need the following properties of the
functional $\tilde{F}_{\omega}$ which establishes the relations
between those functionals.

\bigskip
\noindent {\bf Proposition 1} a).{\it$\tilde{F}_{\omega}(\phi)$
satisfies the cocyle condition:
$\tilde{F}_{\omega}(\phi)+\tilde{F}_{\omega_{\phi}}(\varphi)
=\tilde{F}_{\omega}(\phi+\varphi)$

 $\ \ \ \ \ \ \ \ \ \ \ \ \ \ \
\ \
$b).$\tilde{F}_{\omega}(\phi_{t})=-\frac{1}{t}\int_{0}^{t}(\tilde{I}_{\omega}-\tilde{J}_{\omega})(\phi_{s})ds
-\log(\frac{1}{V}\int_{M}e^{{h}_{\omega}-\phi_{t}}{\omega}^n)$}

\medskip
\noindent Proof of Proposition 1 a). The proof is similar to the
classical case.\\
b). By the definition of $\tilde{F}_{\omega}$, it suffices to
prove that
$\frac{d}{dt}[t(\tilde{J}_{\omega}(\phi_{t})-\frac{1}{V}\int_{M}\phi_{t}e^{-\sigma(u_\omega)}{\omega}^n]
=-(\tilde{I}_{\omega}-\tilde{J}_{\omega})(\phi_{t})$.
By direct computation,
\begin{eqnarray*}
\frac{d}{dt}[t(\tilde{J}_{\omega}(\phi_{t})-\frac{1}{V}\int_{M}\phi_{t}e^{-\sigma(u_\omega)}{\omega}^n)]&=&
\frac{d}{dt}(t\frac{1}{V}\int_{0}^{t}\int_{M}\dot{\phi_{s}}e^{-\sigma(u_{\omega_s})}{\omega}_{\phi_s}^n ds)\\
&=&
\frac{1}{V}\int_{0}^{t}\int_{M}\dot{\phi_{s}}e^{-\sigma(u_{\omega_s})}{\omega}_{\phi_s}^n
ds+t\frac{1}{V}\int_{M}\dot{\phi_{t}}e^{-\sigma(u_{\omega_t})}{\omega}_{\phi_t}^n
\end{eqnarray*}
\begin{eqnarray*}
-(\tilde{I}_{\omega}-\tilde{J}_{\omega})(\phi_{t})&=&
\frac{1}{V}\int_{M}\phi_{t}(e^{-\sigma(u_\omega)}{\omega}^n-e^{-\sigma(u_{\omega_t})}{\omega}_{\phi_t}^n)-
\frac{1}{V}\int_{0}^{t}\int_{M}\dot{\phi_{s}}(e^{-\sigma(u_\omega)}{\omega}^n-e^{-\sigma(u_{\omega_s})}{\omega}_{\phi_s}^n)ds\\
&=&\frac{1}{V}\int_{0}^{t}\int_{M}\dot{\phi_{s}}e^{-\sigma(u_{\omega_s})}{\omega}_{\phi_s}^n
ds-\frac{1}{V}\int_{M}{\phi_{t}}e^{-\sigma(u_{\omega_t})}{\omega}_{\phi_t}^n
\end{eqnarray*}
Differentiating equation (2.4) with respect to t and using the
relation (2.8) one gets
$-\tilde{\Box_{\omega}}_{\phi_t}\dot{\phi}=\phi+t\dot{\phi}$. By
(2.6) we derive that $0=-\int_M
(\tilde{\Box_{\omega}}_{\phi_t}\dot{\phi})e^{-\sigma(u_{\omega_t})}{\omega}_{\phi_t}^n=\int_M
(\phi+t\dot{\phi}) e^{-\sigma(u_{\omega_t})}{\omega}_{\phi_t}^n$.
Hence
$-\frac{1}{V}\int_{M}{\phi_{t}}e^{-\sigma(u_{\omega_t})}{\omega}_{\phi_t}^n=t\frac{1}{V}\int_{M}\dot{\phi_{t}}e^{-\sigma(u_{\omega_t})}{\omega}_{\phi_t}^n$
and the desired equation holds.

\medskip

Now we are in position to prove Theorem 1.

\medskip
\noindent Proof of Theorem 1: By Proposition 1,
$\tilde{F}_{\omega}(\phi_{t})=-\frac{1}{t}\int_{0}^{t}(\tilde{I}_{\omega}-\tilde{J}_{\omega})(\phi_{s}))ds-
\log(\frac{1}{V}\int_{M}e^{{h}_{\omega}-\phi_{t}}{\omega}^n)$.
Also from fact 1 the first term is negative, hence for
$t>\varepsilon$ where $\varepsilon$ is a fixed positive constant,
\begin{eqnarray*} \tilde{F}_{\omega}(\phi_{t})&\leq&
-\log(\frac{1}{V}\int_{M}e^{{h}_{\omega}-\phi_{t}}{\omega}^n)\\
&=&-\log(\frac{1}{V}\int_{M}e^{(t-1)\phi_{t}}{\omega}_{\phi_{t}}^n)\\
&\leq&\frac{1-t}{V}\int_{M}\phi_{t}{\omega}_{\phi_{t}}^n\\
&=&\frac{1-t}{V}\int_{M}\phi_{t}e^{-t\phi_{t}+{h}_{\omega}+\sigma(u_{\omega_t})}{\omega}^n\\
&\leq&\frac{1-t}{V}\int_{ \{ \phi_{t}>
0 \}}\phi_{t}e^{-t\phi_{t}+{h}_{\omega}+\sigma(u_{\omega_t})}{\omega}^n\\
&\leq& C
\end{eqnarray*}
where we used the concavity of  $\log$  and the fact that
$xe^{-tx}$ is uniformly bounded for $x>0$. Here the constant $C$
only depends on the choice of $\varepsilon$ and the initial
metric. Hence by the properness of the
$\tilde{F}_{\omega}(\phi_{t})$ the functional
$\tilde{J}_{\omega}(\phi_{t})$ is uniformly bounded for
$t>\varepsilon$. Consequently we have the bound for $osc(\phi_t)$.
Next, consider the equation (2.4), which implies that
\begin{eqnarray*}
\int_{M}e^{-\sigma(u_{\omega_t})}{\omega}_{\phi_t}^n=V=\int_{M}{\omega}_{\phi_t}^n=\int_{M}e^{{h}_{\omega}-t\phi_t+\sigma(u_{\omega_t})}{{\omega}}^{n}
\end{eqnarray*}
By mean value theorem, there exists a point $x_t$ on $M$ such that
${h}_{\omega}(x_t)-t\phi_t(x_t)+\sigma(u_{\omega_t}(x_t))=0$ for
each time $t$. Hence, we get
\begin{eqnarray*}
|{h}_{\omega}(x)-t\phi_t(x)+\sigma(u_{\omega_t})|&=&|({h}_{\omega}(x)-t\phi_t(x)+\sigma(u_{\omega_t}(x)))-({h}_{\omega}(x_t)-t\phi_t(x_t)+\sigma(u_{\omega_t}(x_t)))|\\
&\leq&tosc(\phi_t)+2|{h}_{\omega}|_{C^0}+2|\sigma|_{C^0} \\
&\leq&C
\end{eqnarray*}
So for $t>\varepsilon$, one gets that $\phi_t$ is uniformly
bounded. The desired $C^0$ estimate is established.

\section{Proof of Theorem 2}
\setcounter{equation}{0}

Assume that there exists an Einstein-Mabuchi metric $\omega_{EM}$
on $M$, In \cite{M4}, Mabuchi has proved that
\medskip

a). The Einstein-Mabuchi metric is unique modulo the action of
$Z^0(X)$.
\medskip

b). There exists a one-parameter family of solutions of (2.4)
$\phi_t \in {\mathcal{H}}_{X}$, $0\leq t \leq 1$ such that
$\omega_{\phi_1}=\omega_{EM}$.
\medskip

 Now fix a $K$-invariant potential $\phi \in
{\mathcal{H}}_{X}$ and set
$\omega=\omega_{EM}+\sqrt{-1}{\partial}\bar{\partial}\phi$.
Consider the complex Monge-Amp\`{e}re equation: \be
\omega_{\phi_{t}}^n=e^{{h}_{\omega}+\sigma(u_{\omega_{\phi_{t}}})-t\phi_{t}}\omega^n
\ee By above result, we have $K$-invariant solution for all
$t\in[0,1]$ and $\omega_{\phi_{1}}=\omega_{EM}$. In particular
$\phi_{1}$ and $-\phi$ differ by a constant.
\medskip

From the section 2 we have seen that the generalized functionals
have similar properties as the classical ones. So by similar
computation as in \cite{PSSW}, we can derive similar inequalities
for $\tilde{F}$, and $\tilde{J}$ which are  \be
|\tilde{J}_{\omega}(\phi_{1})-\tilde{J}_{\omega}(\phi_{0})|\leq
2osc(\phi_{1}-\phi_{0}) \ \ \ \ \ for\ \phi_{0},\ \phi_{1}\in
{\mathcal{H}}_{X}\ee \be \tilde{F}_{\omega_{EM}}(\phi)\geq
c_{0}(1-t)\tilde{J}_{\omega_{EM}}(\phi)-c_{1}(1-t)osc(\phi_{t}-\phi_{1})\ee
where $c_{0},\ c_{1}$ are constant only depending on the choice of
$\sigma$.\\\\
Then we need to estimate $||\phi_{t}-\phi_{1}||_{C^0}$. Rewrite
the equation by using $\omega_{EM}$ as the reference metric, \be
\log{\frac{\omega_{EM}^n}{[\omega_{EM}-\sqrt{-1}{\partial}\bar{\partial}(\phi_{1}-\phi_{t})]^n}}+
(\phi_{1}-\phi_{t})+\sigma(u_{\omega_{\phi_{t}}})-\sigma(u_{\omega_{\phi_{1}}})=(t-1)\phi_{t}\ee
by using the relation
$u_{\omega_{\phi_{t}}}=u_{\omega_{\phi_{1}}}-\sqrt{-1}\bar{X}(\phi_{t}-\phi_{1})$,
we have\be
\log{\frac{\omega_{EM}^n}{[\omega_{EM}-\sqrt{-1}{\partial}\bar{\partial}(\phi_{1}-\phi_{t})]^n}}+
(\phi_{1}-\phi_{t})+\sigma(u_{\omega_{\phi_{1}}}+\sqrt{-1}\bar{X}(\phi_{1}-\phi_{t}))-\sigma(u_{\omega_{\phi_{1}}})=(t-1)\phi_{t}\ee
The linearization of the left hand side of equation at
$\psi=\phi_1-\phi_t=0$ is $\delta\psi$ $\rightarrow
\tilde{\Box}_{\omega_{EM}}\delta\psi+\delta\psi$. Consider the
following space,
\begin{eqnarray*}
\Lambda_1(M, \omega_{EM})=\{u\in C^{\infty}(M) |
\tilde{\Box}_{\omega_{EM}}u=\Delta_{\omega_{EM}}u+\sqrt{-1}\dot{\sigma}(u_{\omega_{EM}})\bar{X}u=-u\}
\end{eqnarray*}
Similar to the classical case, one can prove that $\Lambda_1(M,
\omega_{EM})$ is isomorphic to a subspace of all homomorphic
vector fields on $M$. Assume that  $u_1,\ u_2,\ ...u_m$ form a
basis of this space. Define matrix $k_{ij}(g)$ for each $g\in K$
by $\rho(g)u_i=k_{ij}(g)u_j$, where $\rho$ is the action of $K$ on
space $\Lambda_1(M, \omega_{EM})$. Under our assumption that
$\phi$ is $K$-invariant, we have
\begin{eqnarray*}
V_i=\int_M(\phi
u_i)e^{-\sigma(u_{\omega_{EM}})}{\omega}^n_{EM}=k_{ij}(g)\int_M\phi
u_j e^{-\sigma(u_{\omega_{EM}})}{\omega}^n_{EM}= k_{ij}(g)V_j
\end{eqnarray*}
This implies that the vector $V=(V_1,\ V_2,...,\ V_m)$ is fixed by
$K$. Since $K$ has finite centralizer in $G$ whose Lie algebra is
the set of all holomorphic vector fields on $M$, the vector $V$
must be 0. This is equivalent to say that all $K$-invariant
functions are perpendicular to the space $\Lambda_1(M,
\omega_{EM})$. So the linerized operator is invertible for all
$K$-invariant functions. So we can apply the implicit function
theorem to estimate $||\phi_{t}-\phi_{1}||_{C^0}$ in terms of
$(t-1)\phi_{t}$. Following \cite{PSSW},  we need to prove the the
following which is similar to lemma 1 in \cite{PSSW}, i.e., \be
||\phi_{t}-\phi_{1}||_{C^{0}}\leq
C[(1-t)||\phi_{t}||_{C^{0}}+1]\ee for all $t\in[t_{0},1]$, where
$t_{0}$ (depending on $\phi$) is defined by \be
(1-t_{0})^{1-\alpha}
(1+2(1-t_{0})||\phi_{t_{0}}||_{C^0})^\alpha=\sup_{t\in
[t_{0},1]}(1-t)^{1-\alpha} (1+2(1-t)||\phi_t||_{C^0})^\alpha=D.\ee
and $D$ is a constant only depending on the choice of $p$ and
$\kappa$. Here $p>2n$, $0<\kappa<1$ and
$\alpha=\frac{p+\kappa-2}{p-1}$.

\medskip
\noindent {\bf Remark:} With this bound and Theorem [M1] in the
previous section, one can prove Theorem 2 by using the same
argument in \cite{PSSW} since the properties of the functionals
are similar to the classical case. Since we have an additional
term $\sigma$ in this case, we need to estimate
${h}_{\omega_t}+\sigma(u_{\omega_t})$ and use the volume form
$e^{-\sigma(u_{\omega_t})}{\omega}_t^n$ instead of $h_{\omega_t}$
and $\omega_t^n$. The most important tool to derive the above
bound in \cite{PSSW} is the K\"ahler-Ricci flow. We will introduce
a heat flow and  derive the same smoothing lemma for
${h}_{\omega_t}+\sigma(u_{\omega_t})$ in this case. The bound for
$||\phi_{t}-\phi_{1}||_{C^0}$ is an easy consequence which can be
proved by the method in \cite{PSSW}.

\medskip

For each $t$ consider the following heat flow $f_{s,t}$ in time
$s$ with initial data $f_{0,t}=0$, \be \frac{\partial
f_{s,t}}{\partial s}=\log
\frac{(\omega_{\phi_t}+\sqrt{-1}\partial{\bar{\partial}}f_{s,t})^n}{\omega_{\phi_t}^n}-
{h}_{\omega_{\phi_t}}+f_{s,t}-\sigma(u_{\omega_{\phi_t+f_{s,t}}})\ee
which is the same as \be\frac{\partial
{\omega_{\phi_t+f_{s,t}}}}{\partial
s}=-Ric({\omega}_{\phi_t+f_{s,t}})+\omega_{\phi_t+f_{s,t}}-\sqrt{-1}\ddb{\sigma(u_{\omega_{\phi_t+f_{s,t}}})}
\ee Write $f_t$ for $f_{1,t}$, and consider the K\"ahler form \be
\omega_{\phi_t+f_t}=\omega+
\sqrt{-1}\partial{\bar{\partial}}(\phi_t+f_{t})=\omega_{EM}-\sqrt{-1}
\partial{\bar{\partial}}(\phi_1-\phi_t-f_{t})\ee
There exists a constant so that\be
\log{\frac{\omega_{EM}^n}{[\omega_{EM}-\sqrt{-1}{\partial}\bar{\partial}(\phi_1-\phi_t-f_{t})]^n}}+
(\phi_1-\phi_t-f_{t}-a_t)+\sigma(u_{\omega_{\phi_t+f_{t}}})-\sigma(u_{\omega_{\phi_{1}}})={h}_{\omega_{\phi_t+f_t}}+\sigma(u_{\omega_{\phi_t+f_{t}}})
\ee which can be seen easily by applying $\sqrt{-1}
\partial{\bar{\partial}}$ to both sides.

\medskip

 Notice that this heat flow contains an additional term
involving function $\sigma$ which may cause trouble for
computation. But as long as $\sigma$ is convex, we can still
handle it in the computation.  For convenience, let
$\eta_0=\omega_{\phi_t}$,
$\eta_s=\eta_0+\sqrt{-1}{\partial}\bar{\partial}f$,
${h}_{s}+\sigma(u_s)={h}_{\eta_{s}}+\sigma(u_{\eta_s})$. Then
${h}_{s}+\sigma(u_s)=-\dot{f}+c_s$ for some constant $c_s$ with
$c_0=0$. We will use  s to indicate norms that are defined with
respect to the metric $\eta_s$. Then we prove

\bigskip
\noindent {\bf Lemma 1}{\it We have the following inequalities}:
\begin{eqnarray*}
(a)&\ \ \ \ \ \ \ \ &||\dot{f}||_{C^0}\leq e^s||{h}_0+\sigma(u_0)||_{C^0}\\
(b)&\ \ \ \ \ \ \ \ & \sup_{M}(|\dot{f}|^2+s|\nabla
\dot{f}|_{s}^{2})\leq e^{2s}||{h}_0+\sigma(u_0)||_{C^0}^2\\ (c)&\
\ \ \ \ \ \ \ & \tilde{\Box}_s({h}_s+\sigma(u_s))\geq
e^s\tilde{\Box}_0({h}_0+\sigma(u_0))
\end{eqnarray*}
\\
Proof of Lemma 1.Differentiating the flow one get \be
\frac{\partial}{\partial
s}\dot{f}=\tilde{\Box}_s\dot{f}+\dot{f},\ee hence
$||\dot{f}||_{C^0}\leq e^s||{h}_0+\sigma(u_0)||_{C^0}$, giving (a).
\medskip

Similarly, we compute the flow for $|\nabla \dot{f}|_{s}^{2}$
\begin{eqnarray*}
\frac{\partial}{\partial s}|\nabla
\dot{f}|_{s}^{2}&=&\frac{\partial}{\partial
s}(g^{i\bar{j}}\dot{f}_i\dot{f}_{\bar{j}})\\
&=&-g^{i\bar{a}}(\frac{\partial}{\partial
s}g_{\bar{a}}b)g^{b\bar{j}}\dot{f}_i\dot{f}_{\bar{j}}+g^{i\bar{j}}\ddot{f}_i\dot{f}_{\bar{j}}+g^{i\bar{j}}\dot{f}_i\ddot{f}_{\bar{j}}
\end{eqnarray*}
Use the flow we get that $\frac{\partial}{\partial
s}g_{\bar{a}b}=-R_{\bar{a}b}+g_{\bar{a}b}-\partial_b\partial_{\bar{a}}\sigma(u_{\omega_f})$,
so the first term in the above equation becomes
\begin{eqnarray*}
-g^{i\bar{a}}(\frac{\partial}{\partial
s}g_{\bar{a}b})g^{b\bar{j}}\dot{f}_i\dot{f}_{\bar{j}}&=&
-g^{i\bar{a}}(-R_{\bar{a}b}+g_{\bar{a}b}-\partial_{\bar{a}}\partial_b\sigma(u_{\omega_f}))g^{b\bar{j}}\dot{f_i}\dot{f_{\bar{j}}}\\
&=&-|\nabla
\dot{f}|_{s}^{2}+g^{b\bar{j}}R^{i}_{b}\dot{f}_i\dot{f}_{\bar{j}}-\sqrt{-1}g^{i\bar{a}}(\sqrt{-1}\dot{\sigma}g^{b\bar{j}}\partial_b{u_{\omega_f}})_{\bar{a}}
\dot{f}_i\dot{f}_{\bar{j}}\\
&=&-|\nabla
\dot{f}|_{s}^{2}+g^{b\bar{j}}R^{i}_{b}\dot{f}_i\dot{f}_{\bar{j}}-\sqrt{-1}g^{i\bar{a}}(\dot{\sigma}\overline{X^j})_{\bar{a}}
\dot{f}_i\dot{f}_{\bar{j}}\\
&=&-|\nabla
\dot{f}|_{s}^{2}+g^{b\bar{j}}R^{i}_{b}\dot{f}_i\dot{f}_{\bar{j}}-\sqrt{-1}\ddot{\sigma}g^{i\bar{a}}\partial_{\bar{a}}u_{\omega_{f}}\dot{f}_i(\bar{X}\dot{f})
-\sqrt{-1}\dot{\sigma}g^{i\bar{a}}\overline{X^j}_{\bar{a}}\dot{f}_i\dot{f}_{\bar{j}}\\
&=&-|\nabla
\dot{f}|_{s}^{2}+g^{b\bar{j}}R^{i}_{b}\dot{f}_i\dot{f}_{\bar{j}}+\ddot{\sigma}(X\dot{f})(\bar{X}\dot{f})-\sqrt{-1}\dot{\sigma}g^{i\bar{a}}\overline{X^j}_{\bar{a}}\dot{f}_i\dot{f}_{\bar{j}}
\end{eqnarray*}
Since $\frac{\partial}{\partial
s}\dot{f}=(\Delta_s+\sqrt{-1}\dot{\sigma}\bar{X})\dot{f}+\dot{f}$
we have
\begin{eqnarray*}
\ddot{f}_i&=&
g^{l\bar{k}}\dot{f}_{\bar{k}li}+\dot{f}_i+\sqrt{-1}\ddot{\sigma}\partial_i{u_{\omega_f}}(\bar{X}\dot{f})
+\sqrt{-1}\dot{\sigma} (\bar{X}\dot{f})_i
\end{eqnarray*}
Then
\begin{eqnarray*}
g^{i\bar{j}}\ddot{f_i}\dot{f_{\bar{j}}}&=&
g^{i\bar{j}}g^{l\bar{k}}\dot{f}_{\bar{k}li}\dot{f_{\bar{j}}}+|\nabla
\dot{f}|_{s}^{2}+\sqrt{-1}\ddot{\sigma}g^{i\bar{j}}\partial_i{u_{\omega_f}}(\bar{X}\dot{f})\dot{f_{\bar{j}}}
+\sqrt{-1}\dot{\sigma}g^{i\bar{j}} (\bar{X}\dot{f})_i\dot{f_{\bar{j}}}\\
&=&g^{i\bar{j}}g^{l\bar{k}}\dot{f}_{l\bar{k}i}\dot{f_{\bar{j}}}+|\nabla
\dot{f}|_{s}^{2}+\ddot{\sigma}(\bar{X}\dot{f})^2+\sqrt{-1}\dot{\sigma}g^{i\bar{j}}\overline{X^\alpha}\dot{f}_{\bar{\alpha}i}\dot{f_{\bar{j}}}
\end{eqnarray*}
where we use the fact that $X$ is a holomorphic vector field and
the relation (1.1) in the last line. Similarly
\begin{eqnarray*}
\ddot{f}_{\bar{j}}&=&
g^{l\bar{k}}\dot{f}_{\bar{k}l\bar{j}}+\dot{f}_{\bar{j}}+\sqrt{-1}\ddot{\sigma}\partial_{\bar{j}}{u_{\omega_f}}(\bar{X}\dot{f})
+\sqrt{-1}\dot{\sigma} (\bar{X}\dot{f})_{\bar{j}}\\
&=&g^{l\bar{k}}\dot{f}_{\bar{j}\bar{k}l}-R^{\bar{m}}_{\bar{j}}\dot{f}_{\bar{m}}+\dot{f}_{\bar{j}}+\sqrt{-1}\ddot{\sigma}\partial_{\bar{j}}{u_{\omega_f}}(\bar{X}\dot{f})
+\sqrt{-1}\dot{\sigma}
(\bar{X}\dot{f})_{\bar{j}}\\
 g^{i\bar{j}}\dot{f}_i\ddot{f}_{\bar{j}}&=&
g^{i\bar{j}}g^{l\bar{k}}\dot{f}_{i}\dot{f_{\bar{j}\bar{k}l}}-g^{i\bar{j}}R^{\bar{m}}_{\bar{j}}\dot{f}_{\bar{m}}\dot{f}_{i}+|\nabla
\dot{f}|_{s}^{2}-\ddot{\sigma}(\bar{X}\dot{f})(X\dot{f})\\
&&+\sqrt{-1}\dot{\sigma}g^{i\bar{j}}\overline{X^\alpha}_{\bar{j}}\dot{f}_{\bar{\alpha}}\dot{f_i}
+\sqrt{-1}\dot{\sigma}g^{i\bar{j}}\overline{X^\alpha}\dot{f}_{\bar{\alpha}\bar{j}}\dot{f_i}
\end{eqnarray*}
Combing these terms we get that
\begin{eqnarray*}
\frac{\partial}{\partial s}|\nabla
\dot{f}|_{s}^{2}&=&g^{i\bar{j}}g^{l\bar{k}}\dot{f}_{i\bar{k}l}\dot{f}_{\bar{j}}+g^{i\bar{j}}g^{l\bar{k}}\dot{f}_{i}\dot{f}_{\bar{j}\bar{k}l}+
|\nabla
\dot{f}|_{s}^{2}+\ddot{\sigma}(\bar{X}\dot{f})^2\\&&+\sqrt{-1}\dot{\sigma}g^{i\bar{j}}\overline{X^a}\dot{f}_{\bar{j}\bar{a}}\dot{f}_i
+\sqrt{-1}\dot{\sigma}g^{i\bar{j}}\overline{X^a}\dot{f}_{i\bar{a}}\dot{f}_{\bar{j}}
\end{eqnarray*}
Also we have
\begin{eqnarray*}
\tilde{\Box_s}|\nabla \dot{f}|_{s}^{2}&=&(\Delta_s+\sqrt{-1}\dot{\sigma}\bar{X})(g^{i\bar{j}}\dot{f_i}\dot{f_{\bar{j}}})\\
&=&g^{i\bar{j}}g^{l\bar{k}}\dot{f}_{i\bar{k}l}\dot{f}_{\bar{j}}+g^{i\bar{j}}g^{l\bar{k}}\dot{f}_{i}\dot{f}_{\bar{j}\bar{k}l}+|\nabla\nabla\dot{f}|^2+|\nabla\bar{\nabla}\dot{f}|^2\\
&&+\sqrt{-1}\dot{\sigma}g^{i\bar{j}}\overline{X^a}\dot{f_{i\bar{a}}}\dot{f_{\bar{j}}}+\sqrt{-1}\dot{\sigma}g^{i\bar{j}}\overline{X^a}\dot{f_{i}}\dot{f_{\bar{j}\bar{a}}}
\end{eqnarray*}
Then the flow for $|\nabla \dot{f}|_{s}^{2}$ is
 \be \frac{\partial}{\partial s}|\nabla
\dot{f}|_{s}^{2}=\tilde{\Box}_s|\nabla
\dot{f}|_{s}^{2}-|\nabla\nabla \dot{f}|_{s}^{2}-|\nabla
\bar{\nabla}
\dot{f}|_{s}^{2}+|\nabla\dot{f}|_{s}^{2}+\ddot{\sigma}(\bar{X}\dot{f})^2\ee

Since along the flow f is also invariant under the $X_{{R}}$, then
$(\bar{X}\dot{f})^2=-(X_{I}\dot{f})^2$ where $X_{I}$ is the
imaginary part of holomorphic vector field $X$. Also by the
convexity of $\sigma$ the last term of the above equation is  less
than 0, so \be \frac{\partial}{\partial s}|\nabla
\dot{f}|_{s}^{2}\leq\tilde{\Box}_s|\nabla
\dot{f}|_{s}^{2}-|\nabla\nabla \dot{f}|_{s}^{2}-|\nabla
\bar{\nabla} \dot{f}|_{s}^{2}+|\nabla\dot{f}|_{s}^{2}\ee Next we
compute the flow for $\dot{f}^2$.
\begin{eqnarray*}
\frac{\partial}{\partial s}\dot{f}^2&=&2\ddot{f}\dot{f}\\
&=&2\tilde{\Box_s}\dot{f}\dot{f}+2\dot{f}^2
\end{eqnarray*}
\begin{eqnarray*}
\tilde{\Box_s}\dot{f}^2&=&(\Delta_s+\sqrt{-1}\dot{\sigma}\bar{X})\dot{f}^2\\
&=&2(\Delta_s+\sqrt{-1}\dot{\sigma}\bar{X})\dot{f}\dot{f}+|\nabla\dot{f}|_{s}^{2}
\end{eqnarray*}
so the flow for $\dot{f}^2$ is

 \be \frac{\partial}{\partial
s}\dot{f}^2=\tilde{\Box}_s\dot{f}^2-2|\nabla\dot{f}|_{s}^{2}+2\dot{f}^2\ee
Combing these two flows,\be\frac{\partial}{\partial
s}(|\dot{f}|^2+s|\nabla
\dot{f}|_{s}^{2})\leq\tilde{\Box}_s(|\dot{f}|^2+s|\nabla
\dot{f}|_{s}^{2})+2(|\dot{f}|^2+s|\nabla \dot{f}|_{s}^{2})\ee The
maximum principle implies  \be \sup_{M}(|\dot{f}|^2+s|\nabla
\dot{f}|_{s}^{2})\leq e^{2s}||{h}_0+\sigma(u_0)||_{C^0}^2 \ee
which
proves (b).

\medskip

For $\tilde{\Box}_s\dot{f}$,
\begin{eqnarray*}
\frac{\partial}{\partial
s}(\tilde{\Box}_s\dot{f})&=&\frac{\partial}{\partial
s}(\Delta_s+\sqrt{-1}\dot{\sigma}\bar{X})\dot{f}\\
&=&
(\tilde{\Box}_s(\Delta_s+\sqrt{-1}\dot{\sigma}\bar{X})\ddot{f}_s-g^{i\bar{a}}(\frac{\partial}{\partial
s}g_{\bar{a}b})g^{b\bar{j}}\dot{f}_{\bar{j}i}-\sqrt{-1}\ddot{\sigma}(X\dot{f})(\bar{X}\dot{f})\\
&=&\tilde{\Box}_{s}^2\dot{f}+\tilde{\Box}_s\dot{f}-\ddot{\sigma}(X\dot{f})(\bar{X}\dot{f})-g^{i\bar{a}}g^{b\bar{j}}\dot{f}_{\bar{j}i}\dot{f}_{\bar{a}b}\\
&=&\tilde{\Box}_{s}^2\dot{f}+\tilde{\Box}_s\dot{f}-\ddot{\sigma}(X\dot{f})(\bar{X}\dot{f})
-|\nabla \bar{\nabla}\dot{f}|_{s}^{2}\
\end{eqnarray*}
where we use the flow for $\dot{f}$ and the fact that
$\frac{\partial}{\partial s}g_{\bar{a}b}=\dot{f}_{\bar{a}b}$.
Hence we get  \be \frac{\partial}{\partial
s}\tilde{\Box}_s\dot{f}=
\tilde{\Box}_{s}^2\dot{f}+\tilde{\Box}_s\dot{f}-\ddot{\sigma}(X\dot{f})(\bar{X}\dot{f})
-|\nabla \bar{\nabla}\dot{f}|_{s}^{2}\leq
\tilde{\Box}_{s}^2\dot{f}+\tilde{\Box}_s\dot{f}\ee and c) also
follows from the maximum principle.\\\\
{\bf Lemma 2} {\it Let
$v=({h}_1+\sigma(u_1))-\frac{1}{V}\int_{M}({h}_1+\sigma(u_1))e^{-\sigma(u_1)}{\eta}_{1}^n$,
then for any $p>2n$, there exists constant $C>0$, depending only
on $\omega_{EM}$, $\sigma$ and p so that \be ||v||_{C^0}\leq
C||{h}_0+\sigma(u_0)||_{C^0}^{\frac{p-2}{p-1}}(1-t)^{\frac{1}{p-1}}\ee}
Proof of Lemma 2. Lemma 1 shows that \be ||v||_{C^0}\leq
2e||{h}_0+\sigma(u_0)||_{C^0}\ee Since $v$ is a real-valued
function and $X_R$-invariant, we also have
\begin{eqnarray*}
\int_M|\nabla v|_{1}^2e^{-\sigma(u_1)}{\eta}_{1}^n&=&-\int_M
v(\tilde{\Box}_1 v)e^{-\sigma(u_1)}{\eta}_{1}^n\\
&=&\int_M
(v-\inf_{M}v)(-\tilde{\Box}_1 v)e^{-\sigma(u_1)}{\eta}_{1}^n\\
&\leq& \int_M (v-\inf_{M}v)\sup_{M}(-\tilde{\Box}_1
v)e^{-\sigma(u_1)}{\eta}_{1}^n\\
&\leq& 2V||v||_{C^0}\sup_{M}(-\tilde{\Box}_1 v)
\end{eqnarray*}
Recall that
${h}_0+\sigma(u_0)={h}_{\omega_{\phi_{t}}}+\sigma(u_{\omega_{\phi_t}})$
and thus $Ric({\eta}_0)+\sqrt{-1}\ddb\sigma(u_{\eta_1})>t{\eta}_0$
which implies that $\Delta_{0}({h}_0+\sigma(u_0))\geq -n(1-t)$.
Also one has
${h}_{\omega_{\phi_{t}}}+\sigma(u_{\omega_{\phi_{t}}})=-(1-t)\phi_t+c_t$
then
\begin{eqnarray*}
|\sqrt{-1}\dot{\sigma}\bar{X}({h}_{\omega_{\phi_{t}}}+\sigma(u_{\omega_{\phi_{t}}}))|&\leq&
C(1-t)|X{\phi_{t}}|\\
&=& C(1-t)|u_{\omega_{\phi_t}}-u_{\omega}|\\
&\leq& C(1-t)
\end{eqnarray*}
where we used the fact that $max_{M}u$ and $min_{M}u$ are
holomorphic invariant and constant $C$ depends on $\sigma$. Hence
$\tilde{\Box}_0
{h}_0+\sigma(u_0)=(\Delta_{0}+\sqrt{-1}\dot{\sigma}\bar{X})({h}_0+\sigma(u_0))\geq-(C+n)(1-t)$.
So by Lemma 1,\be -\tilde{\Box}_1 (h_1+\sigma(u_1))\leq
-(C+n)e(1-t) \ee Substituting in the previous inequality gives\be
\int_M|\nabla v|_{1}^2e^{-\sigma(u_1)}{\eta}_{1}^n\leq
2V_0(C+n)e||v||_{C^0}(1-t)\ee Let $p>2n$. Then some constant $C_i$
depending only on $\omega_{EM}$, $\sigma$, $A$ and $p$,
\begin{eqnarray*}
||v||_{C^0}^p&\leq&
C(\int_M|v|^pe^{-\sigma(u_1)}{\eta}_{1}^n+\int_M|\nabla
v|_{1}^pe^{-\sigma(u_1)}{\eta}_{1}^n)\\
&\leq&C_0(||v||_{C^0}^{p-2}\int_M|v|^2e^{-\sigma(u_1)}{\eta}_{1}^n+(e||{h}_0+\sigma(u_0)||_{C^0})^{p-2}\int_M|\nabla
v|_{1}^2e^{-\sigma(u_1)}{\eta}_{1}^n)\\
&\leq& C_1||{h}_0+\sigma(u_0)||_{C^0}^{p-2}\int_M|\nabla
v|_{1}^2e^{-\sigma(u_1)}{\eta}_{1}^n
\end{eqnarray*}
where we have used the Sobolev inequality, the Poincar\'{e}
inequality and applied (b) of Lemma 1. Here the constants in the
Sobolev and Poincar\'{e} inequalities depend only on $\omega_{EM}$
since the metric $\eta_1$ is equivalent to $\omega_{EM}$. Together
with inequality (3.22), this gives \be ||v||_{C^0}^p\leq
C_2(1-t)||{h}_0+\sigma(u_0)||_{C^0}^{p-2}||v||_{C^0}\ee which is
the
inequality to be proved.

\bigskip

With the help of the above two smoothing lemmas we can use the
same argument as in \cite{PSSW} to derive the bound for
$||\phi_{t}-\phi_{1}||_{C^{0}}$, then the Theorem 2 follows at
once by our remark at the beginning of the section.

\section{A holomorphic invariant of Futaki Type}
\setcounter{equation}{0}

As an analogue of the Futaki invariant, we can also define an
invariant of this type  which can be seen as  an obstruction to
the
existence of Einstein-Mabuchi metrics.

\medskip

Let $\eta(M)$ be the complex Lie algebra which consists of all
holomorphic vector fields on $M$. Then we define the functional
associated to multiplier Hermitian structure below,\be
F^{\sigma}_X(V)= \int_M
V({h}_{\omega}+\sigma(u_\omega))e^{-\sigma(u_\omega)}{\omega}^n,\
V\in\eta(M).\ \ and\ \omega\in{\mathcal{K}}_{X}\ee If there exists
an Einstein-Mabuchi metric on $M$, the above functional vanishes.
When $X=0$, the above functional coincides with the Futaki
invariant. The following Theorem shows that the functional is
well-defined and it is a holomorphic
invariant on M.\\\\
{\bf Theorem 3} {\it The functional $F^{\sigma}_X$ is independent
of
the choice of $\omega$ with $\omega\in {\mathcal{K}}_{X}$}.\\
Proof. Let $\omega^{'}$ be another K\"{a}hler form in $C_1(M)$.
Assume that $\omega_{s}=\omega+\sqrt{-1}\ddb\phi(s)$ where
$\phi(s)\in{\mathcal{H}}_{X}$ for $0\leq s\leq 1$ is a path
connecting $\omega$ and $\omega^{'}$ with $\phi(0)=0$ and
$\omega_1=\omega^{'}$. Along the path, we have \be
\frac{d}{ds}(h_s+\sigma(u_s))=
-\dot{\phi}-\tilde{\Box}_s\dot{\phi}\ee

To the homomorphic vector field $V$ one can associate a smooth
complex-valued function $v$ such that
$V^{\alpha}=g^{\alpha\bar{\beta}}\partial_{\bar{\beta}}v$. Note
here we don't require $v$ to be real-valued as in (1.1). For a
real-valued function $f$, one has that
$V(f)=g^{\alpha\bar{\beta}}\partial_{\bar{\beta}}v\partial_{\alpha}f=(\bar{\partial}v,
\bar{\partial}f)_\omega$ .  By (2.6) and $h_s+\sigma(u_s)$ is
real-valued, we compute the derivative of ${F^{\sigma}_X}$,
\begin{eqnarray*}
\frac{d}{ds}{F^{\sigma}_X}(V)&=&\int_M{V(\frac{d}{ds}(h_s+\sigma(u_s)))}e^{-\sigma(u_{\omega_s})}{\omega}_s^n
+\int_M
{V(h_s+\sigma(u_s))}\tilde{\Box}_s\dot{\phi}e^{-\sigma(u_{\omega_s})}{\omega}_s^n\\
&=&\int_M{(\bar{\partial}v,
\bar{\partial}(-\dot{\phi}-\tilde{\Box}_s\dot{\phi}))_{\omega_s}}e^{-\sigma(u_{\omega_s})}{\omega}_s^n+\int_M
{V(h_s+\sigma(u_s))}\tilde{\Box}_s\dot{\phi}e^{-\sigma(u_{\omega_s})}{\omega}_s^n\\
&=&\int_M
\overline{(\bar{\partial}(-\dot{\phi}-\tilde{\Box}_s\dot{\phi}),
\bar{\partial}v})e^{-\sigma(u_{\omega_s})}{\omega}_s^n+\int_M
{V(h_s+\sigma(u_s))}\tilde{\Box}_s\dot{\phi}e^{-\sigma(u_{\omega_s})}{\omega}_s^n\\
&=&\int_M(\dot{\phi}+\tilde{\Box}_s\dot{\phi}){\tilde{\Box}_sv}e^{-\sigma(u_{\omega_s})}{\omega}_s^n+
\int_{M}{V(h_s+\sigma(u_s))}\tilde{\Box}_s\dot{\phi}e^{-\sigma(u_{\omega_s})}{\omega}_s^n\\
&=&\int_{M}\dot{\phi}\tilde{\Box}_s{v}e^{-\sigma(u_{\omega_s})}{\omega}_s^n+\int_M\tilde{\Box}_s\dot{\phi}\{\tilde{\Box}_sv+V(h_s+\sigma(u_s))\}e^{-\sigma(u_{\omega_s})}{\omega}_s^n
\end{eqnarray*}
We need to use integration by part for the first term, notice that
$\phi$ and $\tilde{\Box}_s\phi$ are real-valued, then
\begin{eqnarray*}
\int_{M}\dot{\phi}\tilde{\Box}_s{v}e^{-\sigma(u_{\omega_s})}{\omega}_s^n&=&
\overline{\int_M\dot{\phi}\overline{\tilde{\Box}_s{v}}e^{-\sigma(u_{\omega_s})}{\omega}_s^n}\\
&=&-\int_M\overline{(\bar{\partial}\dot{\phi},\bar{\partial}v)}_{\omega_s}e^{-\sigma(u_{\omega_s})}{\omega}_s^n\\
&=&-\int_M(\bar{\partial}v,
\bar{\partial}\dot{\phi})_{\omega_s}e^{-\sigma(u_{\omega_s})}{\omega}_s^n\\
&=&\int_M v \overline{\tilde{\Box}_s\dot{\phi}}e^{-\sigma(u_{\omega_s})}{\omega}_s^n\\
&=&\int_M v
{\tilde{\Box}_s\dot{\phi}}e^{-\sigma(u_{\omega_s})}{\omega}_s^n
\end{eqnarray*}
hence
\begin{eqnarray*}
\frac{d}{ds}{F^{\sigma}_X}(V)&=&\int_{M}\tilde{\Box}_s\dot{\phi}{v}e^{-\sigma(u_{\omega_s})}{\omega}_s^n
+\int_M\tilde{\Box}_s\dot{\phi}\{\tilde{\Box}_sv+V(h_s+\sigma(u_s))\}e^{-\sigma(u_{\omega_s})}{\omega}_s^n\\
&=&\int_{M}\tilde{\Box}_s\dot{\phi}\{{v}+\tilde{\Box}_sv+V(h_s+\sigma(u_s))\}e^{-\sigma(u_{\omega_s})}{\omega}_s^n\\
&=&-\int_{M}(\bar{\partial} q,
\bar{\partial}\dot{\phi})e^{-\sigma(u_{\omega_s})}{\omega}_s^n
\end{eqnarray*}
where $q={v}+\tilde{\Box}_sv+V(h_s+\sigma(u_s))$. Now we only need
to prove that $q$ is holomorphic. First notice that
\begin{eqnarray*}
V(h_s+\sigma(u_s))&=&V(h_s)+V(\sigma(u_s))\\
&=&V(h_s)+g^{\alpha\bar{\beta}}\partial_{\bar{\beta}}v\partial_{\alpha}(\sigma(u_s))
\end{eqnarray*}
So we can simplify $q$ as
\begin{eqnarray*}
q&=&{v}+\tilde{\Box}_sv+V(h_s)+g^{\alpha\bar{\beta}}\partial_{\bar{\beta}}v\partial_{\alpha}(\sigma(u_s))\\
&=&v+\Delta_{s}v-g^{\alpha\bar{\beta}}\partial_{\bar{\beta}}v\partial_{\alpha}(\sigma(u_s))+V(h_s)+g^{\alpha\bar{\beta}}\partial_{\bar{\beta}}v\partial_{\alpha}(\sigma(u_s))\\
&=&v+\Delta_{s}v+V(h_s)\\
&=&v+g^{i\bar{j}}v_{\bar{j}i}+g^{i\bar{j}}v_{\bar{j}}h_i
\end{eqnarray*}
where we use the definition for $\tilde{\Box}$ in section 2. Then
\begin{eqnarray*}
q_{\bar{l}}&=&
v_{\bar{l}}+g^{i\bar{j}}v_{\bar{j}i\bar{l}}+g^{i\bar{j}}v_{\bar{j}\bar{l}}h_i+g^{i\bar{j}}v_{\bar{j}}h_{\bar{l}i}\\
&=&v_{\bar{l}}+g^{i\bar{j}}v_{\bar{j}\bar{l}i}-g^{i\bar{j}}v_{\bar{j}}R_{\bar{l}i}+g^{i\bar{j}}v_{\bar{j}}(R_{\bar{l}i}-g_{\bar{l}i})\\
&=&0
\end{eqnarray*}
where we use the fact that $h_s$ is the Ricci potential and
$V=g^{i\bar{j}}v_{\bar{j}}$ is a holomorphic vector field, i.e.,
$g^{i\bar{j}}v_{\bar{j}\bar{l}}=0$. Thus we prove that along the
path the derivative of ${F^{\sigma}_X}(V)$ is 0. The theorem
follows as well.\\\\
{\bf Acknowledgements:} The author is grateful to his advisor
D.H.Phong for his advice and constant support. He also would like
to thank Zuoliang Hou for some helpful discussions.

\end{document}